\documentclass[12pt]{amsart}

\usepackage{color,fancyhdr}
\usepackage[all]{xy}

\usepackage{amsmath,amsthm,mathrsfs,amsfonts,
amssymb,times,latexsym,mathabx}

\usepackage[usenames,dvipsnames,svgnames,table]{xcolor}

\DeclareMathAlphabet{\curly}{U}{rsfs}{m}{n}  

\newtheorem{theorem}{Theorem}[section]
\newtheorem{lemma}{Lemma}[section]

\theoremstyle{definition}

\theoremstyle{problem}

\usepackage{xcolor}

\numberwithin{equation}{section}

\makeatletter
\renewcommand{\pmod}[1]{\allowbreak\mkern7mu({\operator@font mod}\,\,#1)}
\makeatother

\newcommand{\be}{\begin{equation}}
\newcommand{\ee}{\end{equation}}

\renewcommand{\b}{\ensuremath{\beta}}

\renewcommand{\le}{\leqslant}

\renewcommand{\ge}{\geqslant}

\begin{document}

\title[Additive completition of thin sets]
{Additive completition of thin sets}

\author{Jin-Hui Fang}
\address{Department of Mathematics, Nanjing University of Information Science $\&$ Technology, Nanjing 210044, PR China}
\email{fangjinhui1114@163.com}

\author{Csaba S\'{a}ndor*}
\address{Institute of Mathematics, Budapest University of Technology and Economics, Egry J\'ozsef utca 1, 1111 Budapest, Hungary; Department of Computer Science and Information Theory, Budapest University of Technology and Economics, M\H{u}egyetem rkp. 3., H-1111 Budapest, Hungary; MTA-BME Lend\"{u}let Arithmetic Combinatorics Research Group, ELKH, M\H{u}egyetem rkp. 3., H-1111 Budapest, Hungary}
\email{csandor@math.bme.hu}
\thanks{* Corresponding author.}
\thanks{The first author is supported by the National Natural Science Foundation of China, Grant No. 12171246 and the Natural Science Foundation of Jiangsu Province, Grant No. BK20211282. The second author is supported by the NKFIH Grants No. K129335.}
\date{\today}%

\begin{abstract}
Two sets $A,B$ of positive integers are called \emph{exact additive complements}, if $A+B$ contains all sufficiently large integers and $A(x)B(x)/x\rightarrow1$. Let $A=\{a_1<a_2<\cdots\}$ be a set of positive integers. Denote $A(x)$ by the counting function of $A$ and $a^*(x)$ by the largest element in $A\bigcap [1,x]$. Following the work of Ruzsa and Chen-Fang, we prove that, for exact additive complements $A,B$ with
$\frac{a_{n+1}}{na_n}\rightarrow\infty$, we have $A(x)B(x)-x\ge \frac{a^*(x)}{A(x)}+o\left(\frac{a^*(x)}{A(x)^2}\right)$ as $x\rightarrow +\infty$.
On the other hand, we also construct exact additive complements $A,B$ with $\frac{a_{n+1}}{na_n}\rightarrow\infty$ such that
$A(x)B(x)-x\le \frac{a^*(x)}{A(x)}+(1+o(1))\left(\frac{a^*(x)}{A(x)^2}\right)$ holds for infinitely many positive integers $x$.
\end{abstract}
\maketitle

\section{\bf Introduction}

Two sets $A,B$ of positive integers are called \emph{additive complements}, if $A+B$ contains all sufficiently large integers. Let $A=\{a_1<a_2<\cdots\}$ be a set of positive integers. Denote  $A(x)$ by the counting function of $A$ and $a^*(x)$ by the largest element in $A\bigcap [1,x]$. If additive complements $A,B$ satisfy
\begin{eqnarray*}
\frac{A(x)B(x)}{x}\rightarrow1,\end{eqnarray*} then we call such $A,B$ \emph{exact additive complements}. In 2001, Ruzsa \cite{Ruzsa1} introduced the following notation which is powerful during the proof of additive complements: let $m>a_1$ be an integer and $k=A(m)$, denote by $L(m)$ the smallest number $l$ for which there are integers $b_1,\cdots,b_l$ such that the numbers $$a_i+b_j,1\le i\le k, 1\le j\le l$$ contain every residue modulo $m$. Obviously, $L(m)\ge m/k$. Ruzsa also proved in \cite{Ruzsa1} that:\vskip2mm

\noindent\textbf{Theorem A.} If
\begin{eqnarray}\label{11}
\frac{a_{n+1}}{na_n}\rightarrow\infty,\end{eqnarray}
then $A$ has an exact complement.\vskip2mm

\noindent\textbf{Theorem B.} Let $A$ be a set satisfying $\frac{A(2x)}{A(x)}\rightarrow 1$. The following are equivalent:\vskip1mm

\noindent (a) $A$ has an exact complement;\vskip1mm
\noindent (b) $A(m)L(m)/m\rightarrow 1$;\vskip1mm
\noindent (c) there is a sequence $m_1<m_2<\cdots$ of positive integers such that $A(m_{i+1})/A(m_i)\rightarrow 1$ and $A(m_i)L(m_i)/m_i\rightarrow 1$.\vskip3mm

In 2017, Ruzsa \cite{Ruzsa2} further considered exact additive complements, that is:\vskip2mm

\noindent\textbf{Theorem C.} For exact additive complements with $\frac{A(2x)}{A(x)}\rightarrow 1$, we have
\begin{eqnarray*}
A(x)B(x)-x\ge (1+o(1))\frac{a^*(x)}{A(x)}\hskip3mm \mbox{as}\hskip3mm x\rightarrow +\infty.\end{eqnarray*}
\vskip2mm

In 2019, Chen and Fang \cite{ChenFang} improved Theorem C by removing the \emph{exact} condition. Furthermore, Chen and Fang also showed in \cite{ChenFang} that the above Theorem C is the best possible.\vskip2mm

\noindent\textbf{Theorem D.}
There exist exact additive complements $A,B$ with $\frac{A(2x)}{A(x)}\rightarrow 1$ such that
\begin{eqnarray*}
A(x)B(x)-x\le (1+o(1))\frac{a^*(x)}{A(x)}\end{eqnarray*}
holds for infinitely many positive integers $x$.
\vskip2mm

In this paper, under the condition \eqref{11} in \cite{Ruzsa1}, we obtain the following result:\vskip2mm

\begin{theorem}\label{thm1}
For exact additive complements $A,B$ with \eqref{11}, we have \begin{eqnarray}\label{12}
A(x)B(x)-x\ge \frac{a^*(x)}{A(x)}+o\left( \frac{a^*(x)}{A(x)^2}\right) \hskip3mm \mbox{as}\hskip3mm x\rightarrow +\infty.\end{eqnarray}
\end{theorem}
\vskip2mm

On the other hand, we also show that $\frac{a^*(x)}{A(x)^2}$ is the best possible.\vskip2mm

\begin{theorem}\label{thm2}
There exist exact additive complements $A,B$ with \eqref{11} such that
\begin{eqnarray}\label{liminf}
\liminf_{x\to \infty }\frac{A(x)B(x)-x-\frac{a^*(x)}{A(x)}}{\frac{a^*(x)}{A(x)^2}}\le 1
\end{eqnarray}
\end{theorem}
\vskip3mm

\section{Proof of Main Results}\label{sec:proofofMT}

Let
$$\sigma (x, n)=|\{ (a, b) : a+b=n, a, b\le x, a\in A, b\in B\} |$$
and
$$\delta (x, n)=|\{ (a, b) : b-a=n, a, b\le x, a\in A, b\in B\} |.$$
\vskip2mm

The idea in the proof of main results is from \cite{ChenFang}-\cite{Ruzsa2}. We will use the following Ruzsa's lemma during the proof of Theorem \ref{thm1}.
\vskip2mm

\begin{lemma}\cite[Lemma 2.1]{Ruzsa2} Let $U$ and $V$ be finite sets of integers and let
$$\sigma (n)=|\{ (u,v) : u\in U, v\in V, u+v=n\} |$$
and
$$\delta (n)= |\{ (u,v) : u\in U, v\in V, v-u=n\} |.$$
Then
$$\sum_{\sigma (n)>1} (\sigma (n)-1) \ge \frac 1{|U|} \sum_{\delta (n)>1}
(\delta (n)-1).$$
\end{lemma}
\vskip2mm

\noindent\textbf{Proof of Theorem \ref{thm1}.} Assume the contrary.
Suppose that \eqref{12} does not hold. Then there exist a positive number $\delta _0(<1)$ and a sequence $x_1<x_2<\dots<x_k<\dots$ such that
\begin{eqnarray}\label{13}
A(x_k)B(x_k)-x_k\le \frac{a^*(x_k)}{A(x_k)}-\delta _0\frac{a^*(x_k)}{A(x_k)^2}.
\end{eqnarray}
We know that \begin{eqnarray*}A(x_k)B(x_k)-x_k&=&\sum_{\substack{a\le x_k, b\le x_k\\
a\in A, b\in B}} 1 -x_k
=\sum_{n=1}^{2x_k} \sigma (x_k, n)-x_k
=\sum_{\substack{n=1\\ \sigma (x_k, n)\ge 1}}^{x_k} (\sigma
(x_k, n)-1)+ \sum_{\substack{n=x_k+1\\ \sigma (x_k, n)\ge
1}}^{2x_k} \sigma (x_k, n) \\
&=& \sum_{\substack{n=1\\ \sigma (x_k,
n)\ge 1}}^{2x_k} (\sigma (x_k, n)-1) + \sum_{\substack{n=x_k+1\\
\sigma (x_k, n)\ge 1} }^{2x_k} 1
 =\sum_{\substack{n=1\\
\sigma (x_k,
n)>1}}^{2x_k} (\sigma (x_k, n)-1) + \sum_{\substack{n=x_k+1\\
\sigma (x_k, n)\ge 1} }^{2x_k} 1.
\end{eqnarray*}
Since $a^*(x_k)\in A$ and $a^*(x_k)+b>x_k$ for all $b\in B$ with
$x_k-a^*(x_k)<b\le x_k$, we have
$$\sum_{\substack{n=x_k+1\\ \sigma (x_k, n)\ge 1} }^{2x_k} 1\ge
B(x_k)-B(x_k-a^*(x_k)).$$ Thus
\begin{equation*}\label{eq2}
A(x_k)B(x_k)-x_k\ge
\sum_{\substack{n\\ \sigma(x_k, n)>1} } (\sigma
(x_k,n)-1)+B(x_k)-B(x_k-a^*(x_k)).\end{equation*}
It infers from Ruzsa's Lemma that
\begin{equation}\label{eq2}
A(x_k)B(x_k)-x_k\ge
\frac{1}{A(x_k)}\sum_{\substack{n\\ \delta(x_k, n)>1} } (\delta
(x_k,n)-1)+B(x_k)-B(x_k-a^*(x_k)).\end{equation}
Let
\begin{eqnarray*}D=\{(a,b) : a\in A, b\in B, a\le b\le x_k-a^*(x_k) \}.\end{eqnarray*}
Then
\begin{equation}\label{eq100}
\sum_{\substack{n\\ \delta(x_k, n)>1} } (\delta
(x_k,n)-1)=\sum_{\substack{n\\ \delta(x_k, n)\ge 1} } (\delta
(x_k,n)-1)\ge |D|-(x_k-a^*(x_k)+1)
\end{equation}
Now we need a lower bound for $|D|$. We consider the following two cases:
\vskip3mm

{\bf Case 1:} $a^*(x_k)>\frac 12x_k$ for infinitely many $k$. By \eqref{11} we know that
\begin{eqnarray*}A\left(\frac{\delta _0}{5}\frac{a^*(x_k)}{A(x_k)}\right)=A(x_k)-1\hskip3mm  \mbox{for all sufficiently large integers}\hskip3mm k.\end{eqnarray*}
Thus, in this case, by Theorem C and $A(x)B(x)/x\rightarrow 1$ we know that
\begin{eqnarray*}
|D|&\ge & \sum_{\substack{\frac{\delta _0}{5}\frac{a^*(x_k)}{A(x_k)}\le b\le x_k-a^*(x_k)\\ b\in B}} A(b)
\ge  A\left(\frac{\delta _0}{5}\frac{a^*(x_k)}{A(x_k)}\right)\left(B(x_k-a^*(x_k))-B\left(\frac{\delta _0}{5}\frac{a^*(x_k)}{A(x_k)}\right)\right)\nonumber\\
&=&(A(x_k)-1)B(x_k-a^*(x_k))-A\left(\frac{\delta _0}{5}\frac{a^*(x_k)}{A(x_k)}\right)B\left(\frac{\delta _0}{5}\frac{a^*(x_k)}{A(x_k)}\right)\nonumber\\
&=& A(x_k)B(x_k)+A(x_k)(B(x_k-a^*(x_k))-B(x_k))-B(x_k-a^*(x_k))-A\left(\frac{\delta _0}{5}\frac{a^*(x_k)}{A(x_k)}\right)B\left(\frac{\delta _0}{5}\frac{a^*(x_k)}{A(x_k)}\right)\nonumber\\
&\ge& x_k+\left(1-\frac{\delta _0}{4}\right)\frac{a^*(x_k)}{A(x_k)}+A(x_k)(B(x_k-a^*(x_k))-B(x_k))-B(a^*(x_k))-\frac{\delta _0}{4}\frac{a^*(x_k)}{A(x_k)}\nonumber\\
&\ge& x_k+\left(1-\frac{\delta _0}{4}\right)\frac{a^*(x_k)}{A(x_k)}+A(x_k)(B(x_k-a^*(x_k))-B(x_k))-\left(1+\frac{\delta _0}{4}\right)\frac{a^*(x_k)}{A(x_k)}-\frac{\delta _0}{4}\frac{a^*(x_k)}{A(x_k)}\nonumber\\
&=& x_k-\frac{3\delta _0}{4}\frac{a^*(x_k)}{A(x_k)}+A(x_k)(B(x_k-a^*(x_k))-B(x_k))
\end{eqnarray*} for sufficiently large $k$. It follows from \eqref{eq2} and \eqref{eq100} that
\begin{eqnarray*}
&&A(x_k)B(x_k)-x_k\\
&\ge& \frac{x_k}{A(x_k)}-\frac{3\delta _0}{4}\frac{a^*(x_k)}{A(x_k)^2}+B(x_k-a^*(x_k))-B(x_k)-\frac{x_k-a^*(x_k)+1}{A(x_k)}
+B(x_k)-B(x_k-a^*(x_k) )\\
&>&\frac{a^*(x_k)}{A(x_k)}-\delta _0\frac{a^*(x_k)}{A(x_k)^2}
\end{eqnarray*}
for sufficiently large $k$, a contradiction with \eqref{13}.
\vskip3mm

{\bf Case 2:} $a^*(x_k)\le \frac{x_k}{2}$ for infinitely many $k$. By \eqref{11} we know that
\begin{eqnarray*}A\left(\frac{\delta_0}{4}\frac{a^*(x_k)}{A(x_k)}\right)=A(x_k)-1\hskip3mm  \mbox{for all sufficiently large integers}\hskip3mm k.\end{eqnarray*}
Thus, in this case, by Theorem C and $A(x)B(x)/x\rightarrow 1$ we know that
\begin{eqnarray*}
|D|&\ge & \sum_{\substack{\frac{\delta _0}{2}\frac{a^*(x_k)}{A(x_k)}<b\le x_k-a^*(x_k)\\ b\in B}} A\left(b-\frac{\delta _0}{4}\frac{a^*(x_k)}{A(x_k)}\right)\nonumber\\&=&\sum_{\substack{\frac{\delta _0}{2}\frac{a^*(x_k)}{A(x_k)}<b<a^*(x_k)+\frac{\delta _0}{4}\frac{a^*(x_k)}{A(x_k)}\\
b\in B}} A\left(b-\frac{\delta _0}{4}\frac{a^*(x_k)}{A(x_k)}\right)+
\sum_{\substack{a^*(x_k)+\frac{\delta _0}{4}\frac{a^*(x_k)}{A(x_k)}\le b\le x_k-a^*(x_k)\\b\in B}} A\left(b-\frac{\delta _0}{4}\frac{a^*(x_k)}{A(x_k)}\right)\nonumber\\
&=&(A(x_k)-1)\left(B\left(a^*(x_k)+\frac{\delta _0}{4}\frac{a^*(x_k)}{A(x_k)}\right)-B\left(\frac{\delta _0}{2}\frac{a^*(x_k)}{A(x_k)}\right)\right)\nonumber \\
&&+A(x_k)\left(B( x_k-a^*(x_k))-B\left(a^*(x_k)+\frac{\delta _0}{4}\frac{a^*(x_k)}{A(x_k)}\right)\right)\nonumber\\
&=&A\left(a^*(x_k)+\frac{\delta _0}{4}\frac{a^*(x_k)}{A(x_k)}\right)B\left(a^*(x_k)+\frac{\delta _0}{4}\frac{a^*(x_k)}{A(x_k)}\right)-B\left(a^*(x_k)+\frac{\delta _0}{4}\frac{a^*(x_k)}{A(x_k)}\right)\nonumber\\
&&-A\left(\frac{\delta _0}{2}\frac{a^*(x_k)}{A(x_k)}\right)B\left(\frac{\delta _0}{2}\frac{a^*(x_k)}{A(x_k)}\right)+A(x_k)B(x_k)+A(x_k)(B(x_k-a^*(x_k))-B(x_k))\nonumber\\
&&-A\left(a^*(x_k)+\frac{\delta _0}{4}\frac{a^*(x_k)}{A(x_k)}\right)B\left(a^*(x_k)+\frac{\delta _0}{4}\frac{a^*(x_k)}{A(x_k)}\right)\nonumber\\
&=&A(x_k)B(x_k)+A(x_k)(B(x_k-a^*(x_k))-B(x_k))-
B\left(a^*(x_k)+\frac{\delta _0}{4}\frac{a^*(x_k)}{A(x_k)}\right)\nonumber\\
&&-A\left(\frac{\delta _0}{2}\frac{a^*(x_k)}{A(x_k)}\right)B\left(\frac{\delta _0}{2}\frac{a^*(x_k)}{A(x_k)}\right)\nonumber\\
&\ge& x_k+\left(1-\frac{\delta _0}{10}\right)\frac{a^*(x_k)}{A(x_k)}+A(x_k)(B(x_k-a^*(x_k))-B(x_k))-\left(1+\frac{\delta _0}{10}\right)\frac{a^*(x_k)+\frac{\delta _0}{4}\frac{a^*(x_k)}{A(x_k)}}{A(x_k)}\nonumber\\
&&-\frac{3\delta _0}{5}\frac{a^*(x_k)}{A(x_k)}\nonumber\\
&\ge&x_k-\frac{9\delta _0}{10}\frac{a^*(x_k)}{A(x_k)}+A(x_k)(B(x_k-a^*(x_k))-B(x_k)),\end{eqnarray*}
for sufficiently large $k$. It follows from \eqref{eq2} and \eqref{eq100} that
\begin{eqnarray*}
&&A(x_k)B(x_k)-x_k\\
&\ge& \frac{x_k}{A(x_k)}-\frac{9\delta _0}{10}\frac{a^*(x_k)}{A(x_k)^2}+B(x_k-a^*(x_k))-B(x_k)-\frac{x_k-a^*(x_k)+1}{A(x_k)}+B(x_k)
-B(x_k-a^*(x_k))\\
&>&\frac{a^*(x_k)}{A(x_k)}-\delta_0\frac{a^*(x_k)}{A(x_k)^2}
\end{eqnarray*}
for sufficiently large $k$, a contradiction with \eqref{13}.
\vskip3mm

This completes the proof of Theorem \ref{thm1}.\hfill$\Box$\\

\noindent\textbf{Proof of Theorem \ref{thm2}.} Let $a_1=1$, $a_2=4$. We will construct the
sequence $a_3,a_4,\dots$ with
\begin{eqnarray}\label{09161}a_{n+1}\gg n^4a_n\end{eqnarray} and there exists a sequence $n_1,n_2,\dots$ such that $a_1,a_2,\dots,a_{n_k}$ form a complete residue system modulo $n_k$ and $n_k|a_{n_k}$. We get such a sequence by greedy algorithm: let $n_1=2$, and if $n_1,n_2,\dots ,n_k$ is already defined, then let $n_{k+1}=a_{n_k}$. Since $a_1,\dots,a_{n_k}$ are distinct residues modulo $a_{n_k}$, we can choose $a_{n_k+1},\dots,a_{n_{k+1}}$ such that $a_{m+1}\gg m^4a_m$ for $m=n_k,\dots n_{k+1}-1$, $a_{n_k}|a_{a_{n_k}}$ and $a_1,\dots ,a_{n_{k+1}}$ are complete residue system modulo $n_{k+1}$. \vskip2mm

For every positive integer $k$, we further take
\begin{eqnarray*}
b_1=n_k,\hskip2mm b_2=2n_k,\hskip2mm \dots,\hskip2mm b_{\frac{a_{n_k}}{n_k}}=\frac{a_{n_k}}{n_k}\cdot n_k.\end{eqnarray*} Then \begin{eqnarray*}a_i+b_j,\hskip2mm 1\le i\le p,\hskip2mm 1\le j\le \frac{a_{n_k}}{n_k}\end{eqnarray*} form a complete residue system modulo $a_{n_k}$. It infers from the definition of $L(a_{n_k})$ that
\begin{eqnarray}\label{09162}L(a_{n_k})=\frac{a_{n_k}}{n_k}.\end{eqnarray}
\vskip2mm

For the set $A=\{a_k\}_{k=1}^{\infty}$ and every positive integer $k$, define $q_k$ as follows:
\begin{eqnarray}\label{5} q_k=\lfloor\frac{a_{k+1}}{k^4a_k}\rfloor,\hskip3mm  \mbox{namely},\hskip3mm
q_k\cdot k^4a_k<a_{k+1}\le (q_k+1)\cdot k^4a_k.\end{eqnarray}
Define the same sets $A,B$ as in \cite[Theorem 3]{Ruzsa1} (replacing $m_k$ by $a_k$ instead). Write $A_k=A\bigcap [1,a_k]$. Take $U_k\subseteq [1,a_k]$ such that $|U_k|=L(a_k)$ and $A_k+U_k$ contains every residue module $a_k$. Let
\begin{eqnarray*}
V_k=U_k+\left\{(q_k-1)a_k,q_ka_k,(q_k+1)a_k,\cdots,
\lfloor \frac{q_{k+1}a_{k+1}}{a_k}\rfloor a_k\right\} \hskip3mm\mbox{and}\hskip3mm B=\bigcup_{k=1}^\infty V_k.\end{eqnarray*}
Let $q_ka_k\le x\le q_{k+1}a_{k+1}$. Since the sequence $\{q_k\}_{k=1}^{\infty}$ defined in \eqref{5} still holds the following property: the sequence $\{q_k\}_{k=1}^{\infty}$ is increasing
to infinite by \eqref{09161} and $A(q_ka_k)\sim A(a_k)$ (In fact, $A(q_ka_k)=k=A(a_k)$ by \eqref{5}). By the same proof as in \cite[Theorem 3]{Ruzsa1}, we know that $A,B$ are additive complements and $A(x)B(x)\sim x$. Thus, the set $A$ with \eqref{09161} has an exact complement $B$. Obviously, such $A$ with \eqref{09161} satisfies \eqref{11}.
\vskip2mm

In the following text we will prove that \eqref{liminf} holds for infinitely many $x_k$. For $x$ with $q_ka_k\le x<(q_{k+1}-1)a_{k+1}$, we have $k\le A(x)\le k+1$ and
\begin{eqnarray}\label{6}
B(x)&\le& \left(\lfloor\frac{x}{a_k}\rfloor-q_k+2\right)L(a_k)
+\sum_{i=2}^k\left(\lfloor\frac{q_ia_i}{a_{i-1}}\rfloor-q_{i-1}+2\right)L(a_{i-1}).
\end{eqnarray}
By Theorem B (b) we know that
\begin{eqnarray*}
L(a_{k-1})\le \frac{2a_{k-1}}{k-1} \hskip4mm\mbox{for large}\hskip4mm k.\end{eqnarray*}
 It infers from \eqref{5} we know that
\begin{eqnarray*}
\sum_{i=2}^k\left(\lfloor\frac{q_ia_i}{a_{i-1}}\rfloor-q_{i-1}+2\right)L(a_{i-1})
\le (k-1)\frac{2q_ka_k}{k-1}=O(q_ka_k)=O\left(\frac{a_{k+1}}{(k+1)^4}\right).\end{eqnarray*}
It is easy to see that for large $k$ we have
$$(q_k-2)L(a_k)\le 2\frac{q_ka_k}{k}=O\left(\frac{a_{k+1}}{(k+1)^5}\right).$$
It follows from \eqref{6} that
\begin{eqnarray}\label{7}
B(x)&\le& \frac{x}{a_k}L(a_k)+O\left(\frac{a_{k+1}}{(k+1)^4}\right).\end{eqnarray}
Choose $x_k=a_{n_k+1}$, where $n_k$ is the index satisfying \eqref{09162}. Then by \eqref{7} we have
\begin{eqnarray*}
&&A(x_k)B(x_k)-x_k-\frac{a^*(x_k)}{A(x_k)}\\
&\le&(n_k+1)\frac{x_k}{n_k}-x_k-\frac{x_k}{n_k+1}
+O\left(\frac{x_k}{(n_k+1)^3}\right)\\
&=&\frac{x_k}{A(x_k)^2}+O\left(\frac{x_k}{A(x_k)^3}\right).
\end{eqnarray*}
\vskip3mm

This completes the proof of Theorem \ref{thm2}.\hfill$\Box$\\


\begin{thebibliography}{99}

\bibitem{ChenFang} Y.G. Chen, J.H. Fang, Additive complements with Narkiewicz's condition. Combinatorica 39 (2019), 813-823.

\bibitem{Ruzsa1} I.Z. Ruzsa, Additive completion of lacunary sequences, Combinatorica 21 (2001), 279-291.

\bibitem{Ruzsa2} I.Z. Ruzsa, Exact additive complements, Quart. J. Math. 68 (2017) 227-235.

\end{thebibliography}
\end{document}